\documentclass[10pt,a4paper]{article} 
\usepackage{color}
\usepackage[utf8]{inputenc}

\usepackage{amsmath, amssymb, theorem, latexsym}
\usepackage[english]{babel}

\usepackage{float}

\usepackage{enumerate}
\usepackage{eufrak}
\usepackage{graphicx}

\def\CA{{\mathcal A}}

\def\CD{{\mathcal D}}

\def\CH{{\mathcal H}}

\def\CP{{\mathcal P}}

\def\CS{{\mathcal S}}

\newcommand{\sca}[2]{\langle#1,#2\rangle}

\newcommand{\beq}{\begin{equation}}
\newcommand{\eeq}{\end{equation}}

\newcommand{\Supp}{\textrm{Supp~}}

\renewcommand{\Re}{{\rm Re\,}}
\renewcommand{\Im}{{\rm Im\,}}
\newtheorem{theorem}{Theorem}[section]
\newtheorem{lemma}[theorem]{Lemma}

\newtheorem{remark}[theorem]{Remark}
\newtheorem{corollary}[theorem]{Corollary}

\makeatletter
\@addtoreset{equation}{section}

\makeatother

\title{Spectral instability for even non-selfadjoint anharmonic oscillators}

\author{Raphaël HENRY\\
\small{Université Paris-Sud}\\
\small{15, Rue Georges Clémenceau, 91400 Orsay, France}\\
\small{raphael.henry@math.u-psud.fr
}}

\begin{document}
\maketitle

\begin{abstract}
We study the instability of the spectrum for a class of non-selfadjoint anharmonic oscillators,
estimating the behavior of the instability indices (\emph{i. e.} the norm of spectral projections) 
associated with the large eigenvalues of these oscillators.
More precisely, we consider the operators
\[
\CA(2k,\theta) = -\frac{d^2}{dx^2}+e^{i\theta}x^{2k}
\]
defined on $L^2(\mathbb{R})$, with $k\geq1$ and $|\theta|<(k+1)\pi/2k$.
We get asymptotic expansions for 
the instability indices, extending the results of \cite{Dav2} and \cite{DavKui}.
\end{abstract}
~\\
\textbf{MSC Ckassification :} 34E20, 34L10, 35P05.

\section{Introduction}
It has been known for several years that the spectrum of a non-selfadjoint operator $\CA$, acting on an
Hilbert space $\CH$, can be very unstable under 
small perturbations of $\CA$. In other words, unlike in the selfadjoint case, the norm of
the resolvent of $\CA$ near the spectrum can blow up much faster than the inverse distance to the spectrum.
Equivalently, the spectrum of its perturbations
$\CA+\varepsilon \mathcal{B}$, with $\varepsilon>0$ and any $\mathcal{B}\in\mathcal{L}(\CH)$, 
$\|\mathcal{B}\|\leq1$, is not necessarily included in the set
$
 \{z\in\mathbb{C} : d(z,\sigma(\CA))\leq\varepsilon\}
$.
\\
Let $\lambda\in\sigma(\CA)$ be an isolated eigenvalue of $\CA$, and let $\Pi_\lambda$ denote the
spectral projection associated with $\lambda$. In order to understand the instability of $\lambda$ (in the above sense),
we define the \emph{instability index} of $\lambda$ as the number
\[
 \kappa(\lambda) = \|\Pi_\lambda\|,
\]
see \cite{Dav2}.\\
Of course $\kappa(\lambda)\geq1$ in any case, and $\kappa(\lambda) = 1$ when $\CA$ is selfadjoint.\\
If $\Pi_\lambda$ has rank $1$, that is, if $\lambda$ is simple in the sense of the algebraic multiplicity, we have 
a convenient expression for $\kappa(\lambda)$, which we shall use extensively in the following:
if $u$ and $u^*$ denote respectively eigenvectors of $\CA$ and $\CA^*$ associated with $\lambda$ and $\bar\lambda$, one can 
easily check \cite{AslDav} that
\beq\label{kapparg1gal}
\kappa(\lambda) = \frac{\|u\|\|u^*\|}{|\sca{u}{u^*}|}.
\eeq
To understand the relation between spectral instability and instability indices, we denote by $\sigma_\varepsilon(\CA)$
the \emph{$\varepsilon$-pseudospectra} of $\CA$, that is the family of sets, indexed by $\varepsilon$,
\[
\sigma_\varepsilon(\CA)  =  \left\{z\in\rho(\CA) : \|(\CA-z)^{-1}\|>\frac{1}{\varepsilon}\right\}\cup\sigma(\CA).
\]
From the perturbative point of view, $\sigma_\varepsilon(\CA)$ can be seen as the union of the perturbed spectra, in the
following sense:
\[
 \sigma_\varepsilon(\CA) =  \bigcup_{
\tiny{\begin{array}{c}\mathcal{B}\in\mathcal{L}(L^2),\\ \|\mathcal{B}\|\leq1
\end{array}}
} \sigma(\CA+\varepsilon\mathcal{B}).
\]
This equivalent formulation follows from a weak version of a theorem due to Roch and Silbermann \cite{RoSi}.\\
Instability indices are closely related to the size of $\varepsilon$-pseudospectra around $\lambda$ (see \cite{AslDav}). 
For instance, if $\CA\in\mathcal{M}_n(\mathbb{C})$ is a diagonalizable matrix with distinct eigenvalues 
$\lambda_1,\dots,\lambda_n$, Embree and Trefethen
show \cite{TrEm} that
the $\varepsilon$-pseudospectra are rather well approximated by disks of radius
$\varepsilon\kappa(\lambda_k)$ around the eigenvalues.
More precisely, 
there exists $\varepsilon_0>0$ such that, for all $\varepsilon\in]0,\varepsilon_0[$,
\beq\label{incluPseudo}
 \bigcup_{\tiny{\lambda_k\in\sigma(\CA)}}D(\lambda_k,\varepsilon\kappa(\lambda_k)+\mathcal{O}(\varepsilon^2))\subset
\sigma_\varepsilon(\CA)\subset
\bigcup_{\tiny{\lambda_k\in\sigma(\CA)}}D(\lambda_k,\varepsilon\kappa(\lambda_k)+\mathcal{O}(\varepsilon^2)).
\eeq
In the case of an infinite dimensional space, the validity of this statement should be investigated, as well as the
dependance on $\lambda_k$ of the $\mathcal{O}(\varepsilon^2)$ terms.\\

In the following, we will consider some \textit{anharmonic oscillators} 
\beq\label{oscanh}
\CA(2k,\theta) = -\frac{d^2}{dx^2}+e^{i\theta}|x|^{2k},
\eeq
where $k\geq1$ and $|\theta|<(k+1)\pi/2k$.
These operators are defined on $L^2(\mathbb{R})$ by considering, first on $\mathcal{C}_0^\infty(\mathbb{R})$,
the associated quadratic form, which is sectorial if 
$|\theta|<(k+1)\pi/2k$, see \cite{Dav2}.
As stated in \cite{Dav2}, its spectrum consists of a sequence of discrete simple eigenvalues, denoted 
in nondecreasing modulus order by $\lambda_n = \lambda_n(2k,\theta)$, $\displaystyle{|\lambda_n|\rightarrow+\infty}$,
and the associated instability indices will be denoted by $\kappa_n(2k,\theta)$.\\
All the spectral projections of $\CA(2k,\theta)$ are of rank $1$ (see Lemma $5$ of \cite{Dav2}),
and if $u_n$ denotes an eigenfunction associated
with $\lambda_n(2k,\theta)$, then formula (\ref{kapparg1gal}) yields
\beq\label{kappan}
 \kappa_n(2k,\theta) = \frac{\int_\mathbb{R}|u_n(x)|^2dx}{\left|\int_\mathbb{R}u_n^2(x)dx\right|},
\eeq
since in this case we have $\CA^*\Gamma = \Gamma\CA$ where $\Gamma$ denotes the complex conjugation, and thus $u_n^* = \bar u_n$.\\

E.-B. Davies showed in \cite{Dav2} that $\kappa_n(2k,\theta)$ grows as $n\rightarrow+\infty$ 
faster than any power of $n$ for any $\theta\neq0$, $|\theta|<(k+1)\pi/2k$.
This statement has been improved in the case $k=1$ of the harmonic oscillator (sometimes
referred as the Davies operator),
since E.-B. Davies and A. Kuijlaars showed \cite{DavKui} that $\kappa_n(2,\theta)$ grows exponentially fast as 
$n\rightarrow+\infty$, with an explicit rate $c(\theta)$ :
\beq\label{ResDavKui}
 \lim\limits_{n\rightarrow+\infty}\frac{1}{n}\log \kappa_n(2,\theta) = c(\theta).
\eeq

The purpose of our work is to prove that this last statement actually holds
for any even anharmonic oscillators $\CA(2k,\theta)$, $k\geq1$.
More precisely, we will improve the estimate by getting asymptotic expansions in powers of $n^{-1}$ as 
$n\rightarrow+\infty$.\\
Let us stress that such a growth of the instability indices implies that the family of eigenfunctions of $\CA(2k,\theta)$
can not possess any of the ``good'' properties usually expected. It does not form a basis, neither in the Hilbert sense nor
in the Riesz or Schauder sense (see \cite{Dav2}, \cite{KreSie}). This excludes any hope of decomposing properly an $L^2$
function along the eigenspaces of the operator.\\
We will prove in section \ref{sSG}, however, that the eigenfunctions form complete sets of the $L^2$ space
(Theorem \ref{total}).\\

Before stating the results of our work, let us specify some notation. Given two functions $f$, $g$ 
and a real sequence $(\alpha_j)_{j\geq0}$, we will write
\beq\label{defEquiv1}
 f(\tau)\underset{\tiny{\tau\rightarrow +\infty}}{\sim}g(\tau)\sum_{j=0}^{+\infty}\alpha_j\tau^{-j}
\eeq
to mean that, for all $N\geq1$,
\beq\label{defEquiv2}
 f(\tau) = g(\tau)\left(\sum_{j=0}^N\alpha_j\tau^{-j}+\mathcal{O}(\tau^{-N-1})\right)
\eeq
as $\tau\rightarrow+\infty$.\\
In the case when $f=f(x,\tau)$, $g=g(x,\tau)$ and $\alpha_j = \alpha_j(x,\tau)$ depend on another variable $x$, we will say that (\ref{defEquiv1})
is uniform with respect to $x$ if the remainder term $\mathcal{O}(\tau^{-N-1})$ in (\ref{defEquiv2}) is uniform with respect to $x$.\\
We define likewise the symbol $\underset{\tiny{\tau\rightarrow0}}{\sim}$.\\

The following theorem was anounced in \cite{Hen}.\\

\begin{theorem}\label{equivAnharm}
 Let $k\in\mathbb{N}^*$ and $\theta$ be such that $0<|\theta|<\frac{(k+1)\pi}{2k}$.
If $\kappa_n(2k,\theta)$ denotes the $n$-th instability index of
$
 \CA(2k,\theta) = -\frac{d^2}{dx^2}+e^{i\theta}x^{2k},
$
 then there exist $K(2k,\theta)>0$ and a real sequence $(C^j(2k,\theta))_{j\geq1}$ such that
\beq\label{dvptAnh}
\kappa_n(2k,\theta)\underset{\tiny{n\rightarrow+\infty}}{\sim}\frac{K(2k,\theta)}{\sqrt{n}}e^{c_k(\theta) n}
\left(1+\sum_{j=1}^{+\infty}C^j(2k,\theta)n^{-j}\right),
\eeq
as $n\rightarrow+\infty$, with
\beq\label{ck}
c_k(\theta) = \frac{2(k+1)\sqrt{\pi}\Gamma\left(\frac{k+1}{2k}\right)\varphi_{\theta,k}(x_{\theta,k})}{\Gamma\left(\frac{1}{2k}\right)} >0,
\eeq
where
\begin{eqnarray*}
 x_{\theta,k} &  = & \left(\frac{\tan (|\theta|/(k+1))}{\sin (k|\theta|/(k+1))+\cos(k|\theta|/(k+1))
\tan (|\theta|/(k+1))}\right)^{\frac{1}{2k}},\\
\varphi_{\theta,k}(x) & = & \Im\int_0^{x e^{i\frac{\theta}{2(k+1)}}}(1-t^{2k})^{1/2}dt~.
\end{eqnarray*}
\end{theorem}
In \cite{Hen}, the well-known asymptotic properties of the Airy function \cite{AbrSteg} have been used to obtain a similar asymptotic expansion for the
instability indices of the complex Airy operator
\[
 -\frac{d^2}{dx^2}+e^{i\theta}|x|,
\]
which can be decomposed as its Dirichlet and Neumann realizations in $\mathbb{R}^+$. It is an example of an odd non-selfadjoint anharmonic oscillator.
The other odd cases $\CA(2k+1,\theta) = -\frac{d^2}{dx^2}+e^{i\theta}|x|^{2k+1}$ are excluded from our work because of the singularity at $x=0$ of 
their potential, and because their eigenfunctions can not be expressed as easily as those of the complex Airy operator in terms of special functions.
However, the instability indices of odd anharmonic oscillators are expected to behave as in (\ref{dvptAnh}).

\begin{remark}
In the harmonic case $k=1$ (Davies operator), we recover from the first term of (\ref{dvptAnh}) the Davies-Kuijlaars theorem
\cite{DavKui} :
\[
 \lim\limits_{n\rightarrow+\infty}\frac{1}{n}\log\|\Pi_n\|  =  c_1(\theta) = 
4\varphi_1\left(\frac{1}{\sqrt{2\cos(\theta/2)}}\right)
 = 2\Re f\left(\frac{e^{i\theta/4}}{\sqrt{2\cos(\theta/2)}}\right)
\]
where
$
 f(z) = \log(z+\sqrt{z^2-1})-z\sqrt{z^2-1}$.
\end{remark}
\vspace{.5cm}
We are also interested in the completeness of the family of eigenfunctions of operator $\CA(2k,\theta)$. The following theorem
has been proved in \cite{Alm} in the case of Airy operator $\CA^D(1,\theta)$, and in \cite{Dav2} in the harmonic oscillator 
case, as 
well as for $\CA(2k,\theta)$, 
$k\geq2$, $|\theta|<\frac{\pi}{2}$. We extend the result to any operator $\CA(2k,\theta)$ with $|\theta|<\frac{(k+1)\pi}{2k}$ :
\begin{theorem}\label{total}
For all $k\geq1$, and for $|\theta|<\frac{(k+1)\pi}{2k}$,
 the eigenfunctions of $\CA(2k,\theta)$
 form a complete set of the space $L^2(\mathbb{R})$.
\end{theorem}

Theorem \ref{total} and the previous estimates enable us 
to study the convergence of the operator series defining the semigroup
$e^{-t\CA(2k,\theta)}$ associated with $\CA(2k,\theta)$
when decomposed along the projections $\Pi_n(2k,\theta)$.\\
The following statement extends the result of \cite{DavKui} in the harmonic oscillator case.

\begin{corollary}\label{seriesemigpes}
Let $|\theta|\leq\pi/2$, $e^{-t\CA(2k,\theta)}$ be the semigroup generated by $\CA(2k,\theta)$,
$\displaystyle{\lambda_n = \lambda_n(2k,\theta)}$ the eigenvalues of $\CA(2k,\theta)$, and $\Pi_n = \Pi_n(2k,\theta)$ the
associated spectral projections.\\
Let $\displaystyle{T(\theta) = c_1(\theta)/\cos(\theta/2)}$, where $c_1(\theta)$ is the constant in (\ref{ck}).
The series
\[
 \Sigma_{2k,\theta}(t) = \sum_{n=1}^{+\infty}e^{-t\lambda_n(2k,\theta)}\Pi_n(2k,\theta)
\]
is not normally convergent in the case $k=1$ for $t<T(\theta)$ ;
in cases $k=1$ for $t>T(\theta)$, and $k\geq2$ for any $t>0$,
 the series converges normally towards $e^{-t\CA(2k,\theta)}$.
\end{corollary}

We prove Theorem \ref{equivAnharm} in Section \ref{sWKB}, while Section \ref{sSG} is dedicated to the proof of Theorem \ref{total} and Corollary
\ref{seriesemigpes}.\\

\textbf{Acknowledgments}\\
I am greatly indebted to Bernard Helffer for his help, advice and comments ; I am also grateful to Thierry Ramond,
Christian Gérard and André Martinez for their valuable discussions. I acknowledge the support of the ANR NOSEVOL.

\section{Instability of even anharmonic oscillators}\label{sWKB}
We would like to understand the behavior as $n\rightarrow+\infty$ of the instability indices $\kappa_n(2k,\theta)$, $k\geq1$,
using Formula (\ref{kappan}). 
In this purpose, we will reformulate the problem in terms of the elements
$\psi_h$ of the kernel of the selfadjoint operator 
\beq\label{defPh}
\CP_h(2k) = -h^2\frac{d^2}{dx^2}+x^{2k}-1.
\eeq

\subsection{Asymptotics of the eigenfunctions}\label{pCAtoCP}
Let $n\geq1$, $\lambda_n$ the $n$-th eigenvalue of $\CA(2k,\theta)$, and $u_n$ an associated eigenfunction.
Let us denote
\beq\label{hn}
h_n = |\lambda_n|^{-\frac{k+1}{2k}}.
\eeq
Since $u_n$ extends to an entire function of the complex plane, we can perform
the analytic dilation $y= e^{i\theta/(2k+2)}|\lambda_n|^{-1/2k}x$, which maps the equation
\[
(\CA(2k,\theta)-\lambda_n)u_n(x) = 0
\]
into
\[
 |\lambda_n|e^{i\theta/(k+1)}\left(-h_n^2\frac{d^2}{dy^2}+y^{2k}-e^{i(\arg\lambda_n-\theta/(k+1))}\right)\psi_{h_n}(y) = 0,
\]
where
\beq\label{defPsiH}
\psi_{h_n}(y) = u_n(h_n^{-1/(k+1)}e^{-i\theta/(2k+2)}y).
\eeq
Since $u_n\in L^2(\mathbb{R})$, according to Sibuya's theory \cite{Sib} (see (\ref{eqwkb}) below), the function $\psi_{h_n}$ is exponentially decreasing in the sectors
$\{|\arg z|<\pi/(2k+2)\}$ and $\{|\arg z-\pi|<\pi/(2k+2)\}$ and it belongs to the domain $H^2(\mathbb{R})\cap L^2(\mathbb{R} ; x^{4k}dx)$. Hence
$\psi_{h_n}$ is an eigenfunction of the non-negative selfadjoint operator
\[
 -h_n^2\frac{d^2}{dy^2}+y^{2k},
\]
associated with the eigenvalue $e^{i(\arg\lambda_n-\theta/(k+1))}$. Therefore, we have necessarily
\[
 e^{i(\arg\lambda_n-\theta/(k+1))} = 1,
\]
which means that $\psi_h$ satisfies the equation
\beq\label{eqPh}
\CP_h(2k)\psi_h = 0,
\eeq
where $\CP_h(2k)$ is the operator defined in (\ref{defPh}).
Furthermore, all the eigenvalues of
$\CA(2k,\theta)$ lie on the half-line $\arg^{-1}\{\frac{\theta}{k+1}\}$.\\
The spectral projection associated with $\lambda_n$ being of rank $1$ according to Lemma $5$ in \cite{Dav2},
Formula (\ref{kappan}) holds for $\kappa_n(2k,\theta)$. Using (\ref{defPsiH}) 
and the scale change $x\mapsto h_n^{\frac{1}{k+1}}x$, and noticing that the solutions $\psi_{h_n}$ are even or odd since 
the potential in $\CA(2k,\theta)$ is even, we get
\beq\label{kappanAnh}
\kappa_n(2k,\theta) = \frac{\int_{\mathbb{R}^+}|\psi_{h_n}(e^{i\frac{\theta}{2(k+1)}}x)|^2dx}
{\int_{\mathbb{R}^+}\psi_{h_n}^2(x)dx}.
\eeq
Here we have deformed the integration path in the denominator, using the analyticity of $\psi_{h_n}$ and its
exponential decay as $|x|\rightarrow+\infty$ in the sector
$\{0\leq\arg z\leq\theta/(2k+2)\}$ \cite{Sib}.
\\
The previous arguments ensure that $r_n = |\lambda_n|$ are the eigenvalues of the selfadjoint anharmonic oscillator
$ -\frac{d^2}{dx^2}+x^{2k}$.
Let us recall from \cite{HelRob}, Theorem $2.1$, the asymptotics of these eigenvalues.
There exists a real sequence $(s_j)_{j\geq1}$ such that
\beq\label{weyl}
|\lambda_n|\underset{\tiny{n\rightarrow+\infty}}{\sim}\left(\frac{(k+1)\sqrt{\pi}\Gamma(\frac{k+1}{2k})}
{\Gamma(\frac{1}{2k})}(n+1/2)\right)^{\frac{2k}{k+1}}\left(1+\sum_{j=1}^{+\infty}s_j(n+1/2)^{-2j}\right).
\eeq
Now we recall some asymptotic properties of the function $\psi_h$.
There exists a canonical domain $\Omega$ for the operator $\CP_h(2k)$, in the sense
of \cite{Sib}, Definition $59.3$, which contains the ray $]0,+\infty[e^{i\theta/(2k+2)}$. Hence, 
according to the results of \cite{Sib} (Theorem $59.1$), since $\psi_h$ is an $L^2$
solution of (\ref{eqPh}), then for any $\delta>0$, we can choose
to normalize $\psi_h$ such that
\beq
\label{eqwkb}
\psi_h(x)\underset{\tiny{h\rightarrow0}}{\sim} \frac{1}{(x^{2k}-1)^{1/4}}
\left(1+\sum_{j=1}^{+\infty}u_j(x)h^j\right)\exp\left(-\frac{1}{h}S(x)\right)
\eeq
uniformly with respect to $x\in[\delta,+\infty[e^{i\theta/(2k+2)}$,
where the functions $u_j$ satisfy, as $|x|\rightarrow+\infty$,
\[
 |u_j(x)| = \mathcal{O}(|x|^{-j(k+1)}).
\]
Here we chose the determination of the square root defined on $\mathbb{C}\setminus[0,+\infty[$, with $\sqrt{-1} = i$,
and we have denoted by $S$ the function
\[
 S : z\mapsto\int_1^z\sqrt{x^{2k}-1}~dx,
\] 
where the integral is taken along a path joining $1$ and $z$, defined on the simply connected set
\[
 D = \mathbb{C}\setminus\bigcup_{j =  0}^{2k-1}e^{ij\pi/k}[1,+\infty[,
\]
the integral being independent of the path.\\

Another expression of $\psi_h$ is available in a complex neighborhood of the real half-line. Namely, there exist $\delta'>0$ and
two sequences  of holomorphic functions $(A_j(\zeta))_{j\geq1}$
and $(B_j(\zeta))_{j\geq0}$ such that
\begin{eqnarray}
 \psi_h(x)& \underset{\tiny{h\rightarrow0}}{\sim} & 2\sqrt{\pi}h^{-1/6}\left(\frac{\zeta(x)}{x^{2k}-1}\right)^{1/4}\left[
 Ai\left(\frac{\zeta(x)}{h^{2/3}}\right)\left(1+\sum_{j=1}^{+\infty}A_j(\zeta(x))h^{2j}\right)\right. \nonumber \\
&& \left.+h^{4/3}Ai'\left(\frac{\zeta(x)}{h^{2/3}}\right)\sum_{j=0}^{+\infty}B_j(\zeta(x))h^{2j}\right]\label{eqAi}
\end{eqnarray}
uniformly with respect to $x\in[-1+\delta',+\infty[+i[-\delta',\delta']$. Here $Ai$ denotes the Airy function and
\[
 \zeta(x) = \left(\frac{3}{2}S(x)\right)^{2/3}.
\]
Notice that $\zeta(x)>0$ if $x>1$, $\zeta(x)\rightarrow+\infty$ as $x\rightarrow+\infty$, and $\zeta(x)<0$ if $x\in]-1,1[$.\\
In order to prove (\ref{eqAi}), let us recall that, according to the results of \cite{Olv}, Theorems $9.1$ and $9.2$, p. $418-419$,
there exists a solution $\tilde\psi_h$ satisfying the asymptotic expansion (\ref{eqAi})
in $\CS_{\delta'} = [-1+\delta',+\infty[+i[-\delta',\delta']$. The strip $\CS_{\delta'}$ is indeed mapped by $\zeta$ into a domain which satisfies 
conditions $(i)-(v)$ of \cite{Olv}, p. $419$. This implies the existence of a solution $\tilde\psi_h$ of (\ref{defPh}) satisfying, for any $n\geq1$,
\[
 \tilde\psi_h(x) = \left(\frac{\zeta(x)}{x^{2k}-1}\right)^{1/4}W_{2n+1,0}(h^{-1},\zeta(x)),
\]
where $W_{2n+1,0}$ is the function given in \cite{Olv}, expression $(9.02)$ p.$~418$. Thus $\tilde\psi_h$ satisfies (\ref{eqAi}).\\
On the other hand, we have
\[
 Ai(\zeta) = \frac{1}{2\sqrt{\pi}\zeta^{1/4}}(1+o(1))\exp\left(-\frac{2}{3}\zeta^{3/2}\right)
\]
as $\zeta\rightarrow+\infty$, see \cite{AbrSteg}.\\
Thus,
\[
 \tilde\psi_h(x) = \frac{1}{(x^{2k}-1)^{1/4}}(1+o(1))\exp\left(-\frac{1}{h}S(x)\right)
\]
as $x\rightarrow+\infty$.
Since $\tilde\psi_h$ is exponentially decaying as $x\rightarrow+\infty$, with the same principal term as $\psi_h$, we have necessarily
$\tilde\psi_h = \psi_h$.\\

\subsection{Estimates on the norm of the eigenfunctions}\label{pEstNorme}
We assume without loss of generality that $\theta>0$ (if $\theta<0$, replace $\theta$ by $|\theta|$).\\
For a fixed $\delta>0$ (which will be determined later in this paragraph),
we write
\beq\label{decompnumer}
\int_0^{+\infty}|\psi_h(e^{i\frac{\theta}{2(k+1)}}x)|^2dx = I_\delta(h)+R_\delta(h)
\eeq
where
\[
I_\delta(h) = \int_{\delta}^{+\infty}|\psi_h(e^{i\frac{\theta}{2(k+1)}}x)|^2dx,~~~
R_\delta(h) = \int_0^{\delta}|\psi_h(e^{i\frac{\theta}{2(k+1)}}x)|^2dx,
\]
and we first estimate $I_\delta(h)$.
The expansion being uniform with respect to $x$, we can take the integral over $[\delta,+\infty[$ in (\ref{eqwkb}).
Thus there exists a sequence $(v_j)_{j\geq1}$ of functions such that
\beq
I_\delta(h) =
\int_{\delta}^{+\infty}a_\theta(x,h)e^{\frac{2}{h}\varphi_{\theta,k}(x)}dx
 \label{DevIdelta}
\eeq
where
\[
 a_\theta(x,h) \underset{\tiny{h\rightarrow0}}{\sim} \frac{1}{|x^{2k}e^{i\frac{k\theta}{k+1}}-1|^{1/2}}
\left(1+\sum_{j=1}^{+\infty}v_j(x)h^j\right)
\]
and
\[
 \varphi_{\theta,k}(x) = -\Re\int_0^{x e^{i\frac{\theta}{2(k+1)}}}(t^{2k}-1)^{1/2}dt.
\]
We have
\[
 \varphi_{\theta,k}'(x) = -|x^{2k}e^{i\frac{k\theta}{k+1}}-1|^{1/2}
\cos\left(\frac{1}{2}\arg(x^{2k}e^{i\frac{k\theta}{k+1}}-1)+\frac{\theta}{2(k+1)}\right)\,.
\]
Hence we can easily check that $\varphi_\theta$ has a unique critical point $x_{\theta,k}$ in $\mathbb{R}^+$,
\beq\label{ptcrit}
x_{\theta,k} = \left(\frac{\tan (\theta/(k+1))}{\sin (k\theta/(k+1))+\cos(k\theta/(k+1))\tan (\theta/(k+1))}\right)^{\frac{1}{2k}}.
\eeq
It is a non-degenerate maximum, and of course $\varphi_{\theta,k}(\delta)<\varphi_{\theta,k}(x_{\theta,k})$ if $\delta<x_{\theta,k}$.\\
Thus, the Laplace method \cite{Erd} applies to the integral (\ref{DevIdelta}), and there exists a sequence $(r_j(2k,\theta))_{j\geq1}$ such that
\beq\label{Idelta}
I_\delta(h) \underset{\tiny{h\rightarrow0}}{\sim} 
\frac{\sqrt{2\pi}}{|(x_{\theta,k}^{2k}e^{ik\theta/(k+1)}-1)\varphi_{\theta,k}''(x_{\theta,k})|^{1/2}}e^{\frac{2}{h}\varphi_{\theta,k}(x_{\theta,k})}h^{1/2}
\left(1+\sum_{j=1}^{+\infty}r_j(2k,\theta)h^j\right).
\eeq

Now the asymptotic expansion (\ref{eqAi}) gives a rough estimate on the remainder term $R_\delta(h)$ in (\ref{decompnumer}), provided that
$\delta$ is chosen small enough. Using the asymptotic behavior of the Airy function and its derivative given in \cite{AbrSteg}
in the sector
$\{|\arg z-\pi|<2\pi/3\}$, expression (\ref{eqAi}) yields, for all $x\in[0,\delta e^{i\theta/(2k+2)}]$,
\[
|\psi_h(x)| = \mathcal{O}(e^{M/h}),~~~M = \sup_{x\in[0,\delta e^{i\theta/(2k+2)}]}(-\Re S(x)).
\]
Choosing $\delta<|x_{\theta,k}|$ then yields
\beq\label{ResteDeltaAnh}
 R_\delta(h) = \mathcal{O}(e^{c/h}),~~~c<2\varphi_{\theta,k}(x_{\theta,k}).
\eeq
Finally (\ref{decompnumer}) and (\ref{Idelta}) lead to the following lemma.
%
\begin{lemma}
There exists a sequence $(r_j(2k,\theta))_{j\geq1}$ such that
\beq\label{numerkappa}
\int_\mathbb{R}|\psi_h(e^{i\frac{\theta}{2(k+1)}}x)|^2dx \underset{\tiny{h\rightarrow0}}{\sim} 
C_k(\theta)h^{1/2}e^{\frac{d_k(\theta)}{h}}\left(1+\sum_{j=1}^{+\infty}r_j(2k,\theta)h^j\right),
\eeq
where
\[
 C_k(\theta) = 2\frac{\sqrt{2\pi}}{|(x_{\theta,k}^{2k}e^{ik\theta/(k+1)}-1)\varphi_{\theta,k}''(x_{\theta,k})|^{1/2}}
 ~~~\textrm{ and }~~~ d_k(\theta) = 2\varphi_{\theta,k}(x_{\theta,k}).
\]
\end{lemma}

In the following paragraph, we get an asymptotic expansion for the denominator of (\ref{kappanAnh}).

\subsection{Estimate on the real axis}
Now we want to get an asymptotic expansion for the norm of $\psi_{h_n}$ on the real axis, using  (\ref{eqAi} which holds in a strip
$[-1+\delta',+\infty[+i[-\delta',\delta']$ for some $\delta'>0$.\\
Let $\chi\in\mathcal{C}^\infty(\mathbb{R};[0,1])$, such that $\Supp\chi\subset]-1+\delta',+\infty[$, and $\chi(-x) = 1-\chi(x)$, $x\in\mathbb{R}$.
Then we have
\[
\int_{\mathbb{R}^+}|\psi_h(x)|^2dx = \int_\mathbb{R}|\psi_h(x)|^2\chi(x)dx.
\]
Noticing that
\[
 Ai\left(\frac{\zeta}{h^{2/3}}\right) = h^{-1/3}Ai_h(\zeta)~~~\textrm{ where }~~~ Ai_h(\zeta) = \int_\mathbb{R}e^{\frac{i}{h}(\zeta\xi+\xi^3/3)}d\xi,
\]
(\ref{eqAi}) yields
\begin{eqnarray}
&&\int_{\mathbb{R}^+}|\psi_h(x)|^2dx  \underset{\tiny{h\rightarrow0}}{\sim} \frac{4\pi}{h}\left(
\int_\mathbb{R}a_1(x,h)|Ai_h(\zeta(x))|^2\chi(x)dx+\right.\nonumber\\
&&+\left.h^4\int_\mathbb{R}a_2(x,h)|Ai_h'(\zeta(x))|^2\chi(x)dx
+h^2\int_\mathbb{R}a_3(x,h)Ai_h(\zeta(x))Ai_h'(\zeta(x))\chi(x)dx
 \right) \nonumber\\
&& =: \frac{4\pi}{h}(I_1(h)+I_2(h)+I_3(h)) \label{decompdenom},
\end{eqnarray}
where for $\ell=1,2,3$,
\[
 a_\ell(x,h)\underset{\tiny{h\rightarrow0}}{\sim}\left|\frac{\zeta(x)}{x^{2k}-1}\right|^{1/2}\sum_{j=0}^{+\infty}a_\ell^j(\zeta(x))h^{2j},~~~
 a_1^0\equiv 1.
\]
In order to estimate $I_1(h)$, we notice (see \cite{Olv}, p. $398$),
that $x\mapsto\zeta(x)$ is one-to-one, mapping $[-1+\delta',+\infty[$ into
$[-\alpha,+\infty[$, for some $\alpha>0$. Let us denote by $x : \zeta\mapsto\zeta(x)$ its inverse, and $\tilde\chi = \chi\circ x$, whose support belongs
to $[-\alpha,+\infty[$. Then,
\begin{eqnarray}
 I_1(h) &  = & \int_\mathbb{R}b_1(\zeta,h)|Ai_h(\zeta)|^2\tilde\chi(\zeta)d\zeta\nonumber \\
  & = & \iiint_{\mathbb{R}^3}e^{\frac{i}{h}\Phi_\xi(\zeta,\eta)}b_1(\zeta,h)\tilde\chi(\zeta)d\zeta,\label{iiint}
\end{eqnarray}
where
\[
 b_1(\zeta,h) \underset{\tiny{h\rightarrow0}}{\sim}\frac{\zeta}{x(\zeta)^{2k}-1}\left(1+\sum_{j=1}^{+\infty}a_1^j(\zeta)h^{2j}\right)
\]
and $\Phi_\xi(\zeta,\eta) = \zeta(\xi-\eta)+(\xi^3-\eta^3)/3$.\\
It is then straightforward to check that the stationary phase method \cite{GriSjo}  applies to
the $(\zeta,\eta)$-integral in (\ref{iiint}), with fixed $\xi$. The unique non-degenerate critical point of $\Phi_\xi$ is
$(\zeta_\xi,\eta_\xi) = (-\xi^2,\xi)$, and we have $\Phi_\xi(\zeta_\xi,\eta_\xi) = 0$, $|\det\textrm{Hess }\Phi_\xi(\zeta_\xi,\eta_\xi)| = 1$. Thus,
there exists a real sequence $(d_j)_{j\geq0}$ such that
\[
 I_1(h)  \underset{\tiny{h\rightarrow0}}{\sim} h\sum_{j=0}^{+\infty}d_jh^j,~~~~d_0>0.
\]
The same treatment for the terms $I_2(h)$ and $I_3(h)$ in (\ref{decompdenom}), using that
\[
 Ai_h'(\zeta) = \frac{i}{h}\int_\mathbb{R}\xi e^{\frac{i}{h}(\zeta\xi+\xi^3/3)}d\xi,
\]
yields
\beq\label{denomkappa}
 \int_{-\infty}^{+\infty}|\psi_h(x)|^2dx 
\underset{\tiny{h\rightarrow0}}{\sim} \sum_{j=0}^{+\infty}c_jh^j,~~~c_0\neq0.
\eeq
\subsection{Proof of Theorem \ref{equivAnharm}}
Finally, we get the desired statement by quantification of the parameter $h_n$ as an asymptotic expansion in powers of $n^{-1}$. Namely,
using (\ref{hn}) and (\ref{weyl}), we have
\[
 \frac{1}{h_n} \underset{\tiny{n\rightarrow+\infty}}{\sim}
\left(\frac{(k+1)\sqrt{\pi}\Gamma(\frac{k+1}{2k})}{\Gamma(\frac{1}{2k})}(n+1/2)\right)
\left(1+\sum_{j=1}^{+\infty}s_k^j(n+1/2)^{-2j}\right).
\]
for some real sequence $(s_k^j)_{j\geq1}$.\\
This expansion along with expressions (\ref{kappanAnh}), (\ref{numerkappa}) and (\ref{denomkappa}) 
yield the statement of Theorem \ref{equivAnharm}.\\

In the last section, we prove Theorem \ref{total} and Corollary \ref{seriesemigpes}.
\section{Completeness and semigroups}\label{sSG}
\subsection{Completeness of eigenfunctions}
In this paragraph we prove Theorem \ref{total}.
First of all, let us recall that, if  $\CH$ is an Hilbert space and $p\geq1$,
the \emph{Schatten class} $C^p(\CH)$ denotes the set of compact operators $\CA$ such that
\beq\label{normSchatten}
 \|\CA\|_p := \left(\sum_{n=1}^{+\infty}\mu_n(\CA)^p\right)^{1/p} <+\infty,
\eeq
where $(\mu_n(\CA))_{n\geq1}$ are the eigenvalues of $(\CA^*\CA)^{1/2}$, repeated according to their multiplicity
(see \cite{DunSch}).
The space $C^p(\CH)$, $p\geq1$, is a Banach space.\\
We already know that the resolvent $\CA(2k,\theta)^{-1}$ is compact for any $k\geq1$ and $|\theta|<(k+1)\pi/2k$.
We now prove like in \cite{PTLR} that it actually belongs to a Schatten class:
\begin{lemma}\label{anharmSchatten}
For any $\varepsilon>0$, $|\theta|<\frac{(k+1)\pi}{2k}$ and $k\geq1$, we have
\[
 (\CA(2k,\theta))^{-1}\in C^{\frac{k+1}{2k}+\varepsilon}(L^2(\mathbb{R})).
\]
\end{lemma}
\textbf{Proof: }
Let us show that, for all $\varepsilon>0$, the series $\sum\mu_n^{\frac{k+1}{2k}+\varepsilon}$ is convergent, where
$(\mu_n)_{n\geq1}$ are the eigenvalues of
\[
 \left([(\CA(2k,\theta))^{-1}]^*(\CA(2k,\theta))^{-1}\right)^{1/2} = \left([\CA(2k,\theta)(\CA(2k,\theta))^*]^{-1}\right)^{1/2}.
\]
If $(\nu_n)_{n\geq1}$ denote the eigenvalues of $\CA(2k,\theta)(\CA(2k,\theta))^*$, then we have to check that
\[
 \sum_{n=1}^{+\infty}\nu_n^{-p/2}<+\infty
\]
as soon as $p>\frac{k+1}{2k}$.\\
$\CA(2k,\theta)(\CA(2k,\theta))^*$ is a selfadjoint operator, and if $p(x,\xi)$ denotes its symbol, we define its
quasi-homogeneous principal symbol $P(x,\xi)$ as
\[
 P(x,\xi) = \lim\limits_{r\to+\infty}r^{-1}p(r^{1/4k}x,r^{1/4}\xi),
\]
following \cite{Rob}.\\
Then we have
\begin{eqnarray}
&& P(x,\xi) = |\xi^2+e^{i\theta}x^{2k}|^2 = \xi^4+2\cos\theta\xi^2x^{2k}+x^{4k},\nonumber \\
 & & P(r^{1/4k}x,r^{1/4}\xi) = rP(x,\xi),~~r>0 .\label{quasiHom}
\end{eqnarray}
Moreover $P$ is globally elliptic, in the sense that
\beq
\forall (x,\xi)\neq(0,0),~~|P(x,\xi)|>0  \label{ellip}.
\eeq
Hence the results of \cite{Rob}, Theorem $7.1$, allow us to apply the following Weyl formula: 
\[
 N(t) := \#\{j\geq1 : \nu_j\leq t\}\underset{\tiny{t\rightarrow+\infty}}{\sim}\int_{P(x,\xi)\leq t}dxd\xi,
\]
which, with $t = \nu_n$ and using (\ref{quasiHom}), yields
\[
 n\underset{\tiny{n\rightarrow+\infty}}{\sim}C\nu_n^{\frac{k+1}{4k}}
\]
where $C=\textrm{Vol }P^{-1}([0,1])$. \\
Thus the series $\sum\nu_n^{-p/2}$ converges if and only if
\[
 \sum_{n=1}^{+\infty}n^{-\frac{2kp}{k+1}}<+\infty,
\]
that is if and only if $p>\frac{k+1}{2k}$.
\hfill $\boxminus$ \\

Since the operator $\CA(2k,\theta)$ is sectorial and its numerical range is included in the sector
$\CS_\theta = \arg^{-1}[0,\theta]$, the resolvent estimate
\beq\label{estRes}
 \|(\CA(2k,\theta)-\lambda)^{-1}\| = \mathcal{O}(|\lambda|^{-1})
\eeq
holds outside $\CS_\theta$, and if we denote $p=\frac{k+1}{2k}+\varepsilon$, then
\beq\label{thetapip}
 \theta < \frac{(k+1)\pi}{2k} < \frac{\pi}{\frac{k+1}{2k}+\varepsilon} = \frac{\pi}{p}
\eeq
for $\varepsilon$ small enough, as soon as $k>1$.\\
Consequently, Theorem \ref{total} follows from Lemma \ref{anharmSchatten} and Corollary
$31$ of \cite{DunSch}, p. 1115.\\

In the next paragraph, we prove Corollary \ref{seriesemigpes}.

\subsection{Semigroup decomposition}
The case $k=1$ was already proved in \cite{DavKui}. For $k\geq2$, using (\ref{dvptAnh}) and (\ref{weyl}), we see that, as $n\rightarrow+\infty$,
\[
 \|\Pi_n(2k,\theta)\| = \mathcal{O}(e^{c|\lambda_n|^\alpha}),
\]
where $\alpha<1$.
Thus, the series $\Sigma_{2k}(t)$ is normally convergent for all $t>0$.\\
To check that the series $\Sigma_{2k}(t)$ (when convergent) converges towards the semigroup associated with $\CA(2k,\theta)$,
we use the density of the family $(u_n)$, where the eigenfunctions
$u_n$ are assumed to be normalized by the condition $\sca{u_n}{\bar u_n} = 1$, so that
$(u_n,\bar u_n)_{n\geq1}$ is a biorthogonal family (see \cite{Dav2}), namely
\beq\label{biortho}
 \forall n,m\in\mathbb{N},~~\sca{u_n}{\bar u_m} = \delta_{n,m}.
\eeq
Then we have
\[
 e^{-t\CA(2k,\theta)}u_n = e^{-t\lambda_n}u_n
\]
and on the other hand,
\[
 \Sigma_{2k}(t)u_n = \sum_{j = 1}^{+\infty}e^{-t\lambda_j}\Pi_ju_n = e^{-t\lambda_n}u_n.
\]
Here we used the formula
\[
\Pi_j f = \sca{f}{\bar u_j}u_j
\]
(see \cite{Dav2}, \cite{AslDav}) which holds for rank $1$ spectral projections,
together with the biorthogonal property (\ref{biortho}).\\
Hence by linearity, $e^{-t\CA(2k,\theta)}$ and $\Sigma_{2k}(t)$ coincide on $\textrm{Vect}\{u_n : n\geq1\}$, 
and hence on $\CD(\CA(2k,\theta))$ by density (see Theorem \ref{total}).\\

\end{document}